\newif\ifsect\newif\iffinal
\secttrue\finaltrue
\def\strutdepth{\dp\strutbox}

\def\lsimb#1{\vadjust{\vtop to0pt{\baselineskip\strutdepth\vss
	\llap{\ttt\string #1\ }\null}}}
\def\rsimb#1{\vadjust{\vtop to0pt{\baselineskip\strutdepth\vss
	\line{\kern\hsize\rlap{\ttt\ \string #1}}\null}}}
\def\ssect #1. {\bigbreak\indent{\bf #1.}\enspace\message{#1}}
\def\smallsect #1. #2\par{\bigbreak\noindent{\bf #1.}\enspace{\bf #2}\par
	\global\parano=#1\global\eqnumbo=1\global\thmno=1
	\nobreak\smallskip\nobreak\noindent\message{#2}}
\def\thm #1: #2{\medbreak\noindent{\bf #1:}\if(#2\thmp\else\thmn#2\fi}
\def\thmp #1) { (#1)\thmn{}}
\def\thmn#1#2\par{\enspace{\sl #1#2}\par
        \ifdim\lastskip<\medskipamount \removelastskip\penalty 55\medskip\fi}

\def\qedn{\thinspace\null\nobreak\hfill\hbox{\vbox{\kern-.2pt\hrule height.2pt 
depth.2pt\kern-.2pt\kern-.2pt \hbox to2.5mm{\kern-.2pt\vrule width.4pt
\kern-.2pt\raise2.5mm\vbox to.2pt{}\lower0pt\vtop to.2pt{}\hfil\kern-.2pt
\vrule width.4pt\kern-.2pt}\kern-.2pt\kern-.2pt\hrule height.2pt depth.2pt
\kern-.2pt}}\par\medbreak}
\def\pf{\ifdim\lastskip<\smallskipamount \removelastskip\smallskip\fi
        \noindent{\sl Proof\/}:\enspace}

\def\forclose#1{\hfil\llap{$#1$}\hfilneg}
\def\newforclose#1{
	\ifsect\xdef #1{(\number\parano.\number\eqnumbo)}\else
	\xdef #1{(\number\eqnumbo)}\fi
	\hfil\llap{$#1$}\hfilneg
	\global \advance \eqnumbo by 1
	\iffinal\else\rsimb#1\fi}
\def\newforclosea#1{
	\ifsect\xdef #1{{\rm(\the\Apptok.\number\eqnumbo)}}\else
	\xdef #1{(\number\eqnumbo)}\fi
	\hfil\llap{$#1$}\hfilneg
	\global \advance \eqnumbo by 1
	\iffinal\else\rsimb#1\fi}
\def\forevery#1#2$${\displaylines{\let\eqno=\forclose\let\neweqa=\newforclosea
        \let\neweq=\newforclose\hfilneg\rlap{$\qqquad\forall#1$}\hfil#2\cr}$$}
\newcount\parano
\newcount\eqnumbo
\newcount\thmno
\newcount\versiono
\versiono=0
\def\neweqt#1$${\xdef #1{(\number\parano.\number\eqnumbo)}
	\eqno #1$$
	\iffinal\else\rsimb#1\fi
	\global \advance \eqnumbo by 1}
\def\newthmt#1 #2: #3{\xdef #2{\number\parano.\number\thmno}
	\global \advance \thmno by 1
	\medbreak\noindent
	\iffinal\else\lsimb#2\fi
	{\bf #1 #2:}\if(#3\thmp\else\thmn#3\fi}
\def\neweqf#1$${\xdef #1{(\number\eqnumbo)}
	\eqno #1$$
	\iffinal\else\rlap{$\smash{\hbox{\hfilneg\string#1\hfilneg}}$}\fi
	\global \advance \eqnumbo by 1}
\def\newthmf#1 #2: #3{\xdef #2{\number\thmno}
	\global \advance \thmno by 1
	\medbreak\noindent
	\iffinal\else\llap{$\smash{\hbox{\hfilneg\string#1\hfilneg}}$}\fi
	{\bf #1 #2:}\if(#3\thmp\else\thmn#3\fi}
\def\inizia{\ifsect\let\neweq=\neweqt\else\let\neweq=\neweqf\fi
\ifsect\let\newthm=\newthmt\else\let\newthm=\newthmf\fi}
\def\bititolo{\empty}
\gdef\begin #1 #2\par{\xdef\titolo{#2}
\ifsect\let\neweq=\neweqt\else\let\neweq=\neweqf\fi
\ifsect\let\newthm=\newthmt\else\let\newthm=\newthmf\fi
\centerline{\titlefont\titolo}
\if\bititolo\empty\else\medskip\centerline{\titlefont\bititolo}
\xdef\titolo{\titolo\ \bititolo}\fi
\bigskip
\centerline{\bigfont \autore}
\if\istituto!\else\bigskip
\centerline{\istituto}
\centerline{\indirizzo}
\centerline{\email}\fi
\medskip
\centerline{#1~\anno}
\bigskip\bigskip
\ifsect\else\global\thmno=1\global\eqnumbo=1\fi}
\def\istituto{Dipartimento di Matematica, Universit\`a di Roma ``Tor Vergata''}
\def\indirizzo{Via della Ricerca Scientifica, 00133 Roma, Italy}
\def\email{E-mail: abate@mat.uniroma2.it}
\def\anno{\the\year}
\def\copertina #1{\vfill\eject\headline={\hfill}
\vbox{\hrule\hbox{\vrule\advance\vsize by-0.8pt\vbox to\vsize{\vfil
\centerline{\bigbigfont Marco Abate}
\vfil
\centerline{\bigsl\titolo}
\vfil\vfil
\centerline{\sns L}
\bigskip
\centerline{\bigbigfont Scuola Normale Superiore}
\medskip
\centerline{\bigbigfont #1~\anno}
\vfil}\vrule}\hrule}}
\font\titlefont=cmssbx10 scaled \magstep1
\font\bigfont=cmr12

\font\bigsl=cmsl12

\font\bigbigfont=cmr10 scaled \magstep2
\font\ttt=cmtt10 at 10truept
\font\eightrm=cmr8

\let\sc=\smallcaps

\font\bbr=msbm10
\font\sbbr=msbm7 
\font\ssbbr=msbm5
\def\ca #1{{\cal #1}}
\nopagenumbers
\binoppenalty=10000
\relpenalty=10000
\newfam\amsfam
\textfont\amsfam=\bbr \scriptfont\amsfam=\sbbr \scriptscriptfont\amsfam=\ssbbr
\let\de=\partial

\def\Im{\mathop{\rm Im}\nolimits}

\def\id{\mathop{\rm id}\nolimits}

\mathchardef\void="083F

\def\C{{\mathchoice{\hbox{\bbr C}}{\hbox{\bbr C}}{\hbox{\sbbr C}}
{\hbox{\sbbr C}}}}
\def\N{{\mathchoice{\hbox{\bbr N}}{\hbox{\bbr N}}{\hbox{\sbbr N}}
{\hbox{\sbbr N}}}}

\def\Q{{\mathchoice{\hbox{\bbr Q}}{\hbox{\bbr Q}}{\hbox{\sbbr Q}}
{\hbox{\sbbr Q}}}}

\def\qqquad{\quad\qquad}

\newcount\notitle
\notitle=1
\headline={\ifodd\pageno\rhead\else\lhead\fi}
\def\rhead{\ifnum\pageno=\notitle\iffinal\hfill\else\hfill\tt Version 
\the\versiono; \the\day/\the\month/\the\year\fi\else\hfill\eightrm\titolo\hfill
\folio\fi}
\def\lhead{\ifnum\pageno=\notitle\hfill\else\eightrm\folio\hfill\autore\hfill
\fi}
\def\autore{Marco Abate}
\output={\plainoutput}
\newbox\bibliobox
\def\setref #1{\setbox\bibliobox=\hbox{[#1]\enspace}
	\parindent=\wd\bibliobox}
\def\biblap#1{\noindent\hang\rlap{[#1]\enspace}\indent\ignorespaces}
\def\art#1 #2: #3! #4! #5 #6 #7-#8 \par{\biblap{#1}#2: {\sl #3\/}.
	#4 {\bf #5} (#6)\if.#7\else, \hbox{#7--#8}\fi.\par\smallskip}
\def\book#1 #2: #3! #4 \par{\biblap{#1}#2: {\bf #3.} #4.\par\smallskip}
\def\coll#1 #2: #3! #4! #5 \par{\biblap{#1}#2: {\sl #3\/}. In {\bf #4,} 
#5.\par\smallskip}
\def\pre#1 #2: #3! #4! #5 \par{\biblap{#1}#2: {\sl #3\/}. #4, #5.\par\smallskip}
\def\End{\mathop{\rm End}\nolimits}
\newcount\defno
\def\smallsect #1. #2\par{\bigbreak\noindent{\bf #1.}\enspace{\bf #2}\par
    \global\parano=#1\global\eqnumbo=1\global\thmno=1\global\defno=0
    \nobreak\smallskip\nobreak\noindent\message{#2}}
\def\newdef #1\par{\global \advance \defno by 1
    \medbreak
{\bf Definition \the\parano.\the\defno:}\enspace #1\par
\ifdim\lastskip<\medskipamount \removelastskip\penalty 55\medskip\fi}
\def\autore{Marco Abate${}^1{}^*$ and Francesca Tovena${}^2$%
\footnote{${}^*$}{\rm The authors have been partially supported by {\it
Centro di Ricerca Matematica ``Ennio de Giorgi",} Pisa.}}
\def\indirizzo{\vbox{\hfill${}^2$Dipartimento di Matematica, Universit\`a
di Roma Tor Vergata,\hfill\break\null\hfill Via della Ricerca
Scientifica, 00133 Roma, Italy.
E-mail: tovena@mat.uniroma2.it\hfill\null}}
\def\istituto{\vbox{\hfill${}^1$Dipartimento di Matematica, Universit\`a
di Pisa,\hfill\break\null\hfill Via Buonarroti 2, 56127 Pisa,
Italy. E-mail: abate@dm.unipi.it\hfill\null\vskip5pt}}
\def\email{}
\begin {August} Formal classification of holomorphic maps tangent
to the identity

{{\narrower{\sc Abstract.}  We describe a procedure for
constructing formal normal forms of holomorphic maps with a hypersurface
of fixed points, and we apply it to obtain a complete list of formal
normal forms for 2-dimensional holomorphic maps tangential to a curve of
fixed points.
\par
}}
\smallsect 0. Introduction

When studying a class of holomorphic dynamical systems, one of the main
goals often is the classification under topological, holomorphic or
formal conjugation; one would like to have a complete list of all the
possible topological (respectively, holomorphic or formal) normal forms.

For discrete holomorphic local dynamical systems tangent to the
identity,  that is germs~$f$ about the origin of holomorphic self-maps
of $\C^n$ such  that $f(O)=O$ and $df_O=\id$, this problem has been
extensively studied  in dimension one.

The formal classification (where the conjugating map is just a formal
power series, not necessarily convergent) is elementary, and
depends on one discrete parameter and one continuous parameter:

\newthm Proposition \formaltangent: Let
$f$ be a $1$-dimensional discrete holomorphic local dynamical system
tangent to the identity of the form 
$$
f(z)=z+a_{\nu+1} z^{\nu+1}+O(z^{r+2}),
$$
where $a_{\nu+1}\ne 0$. Then $f$ is formally conjugated to the map
$$
z\mapsto z+z^{\nu+1}+\beta z^{2\nu+1},
\neweq\eqzuno
$$
where $\beta\in\C$ is a formal (and holomorphic) invariant given
by
$$
\beta={1\over 2\pi i}\int_\gamma{dz\over z-f(z)},
$$
where the integral is taken over a small positive loop~$\gamma$ about the
origin.

The topological classification is even easier to state (but much more
difficult to prove), since it depends only on one discrete parameter:   

\newthm Theorem \Camacho: (Camacho, 1978 [C]; Shcherbakov, 1982 [S]) Let
$f$ be a $1$-dimensional discrete holomorphic local dynamical system
tangent to the identity of the form 
$$
f(z)=z+a_{\nu+1} z^{\nu+1}+O(z^{\nu+2}),
$$
where $a_{\nu+1}\ne 0$. Then $f$ is topologically locally conjugated to
the map
$$
z\mapsto z+z^{\nu+1}.
$$

The holomorphic classification is much more complicated; in particular,
for every choice of $\nu\ge 1$ and~$\beta\in\C$ there are an uncountable
number of $1$-dimensional discrete holomorphic local dynamical system
tangent to the identity which are {\it not} holomorphically conjugated
but have the same formal normal form~\eqzuno. Anyway, \'Ecalle [\'E1--2]
and Voronin [V] have found a complete set of invariants for the
holomorphic classification; see also~[I] and
[M1--2] for a sketch of the proof.

When $n\ge 2$, as far as we know almost nothing is known for the
topological and holomorphic classification. Our note is a contribution
to the formal classification, at least in dimension~2. To explain our
approach, let us recall what is already known. 

In his monumental work [\'E3] (see also [\'E4] for a short survey)
\'Ecalle studied the formal classifications of discrete holomorphic local
dynamical systems tangent to the identity in dimension~$n\ge 2$, giving
a complete set of formal invariants for maps satisfying a generic
condition. To be more precise, write a discrete holomorphic local
dynamical systems tangent to the identity~$f$ in the form
$$
f(z)=z+P_\nu(z)+O(\|z\|^{\nu+1}),
\neweq\eqzumez
$$
where $P_\nu\not\equiv O$ is an $n$-uple of homogeneous polynomials of
total degree~$\nu$ in the variables $z=(z_1,\ldots,z_n)$. A {\sl
characteristic direction} (or {\sl eigenradius,} in \'Ecalle's
terminology) for~$f$ is a direction $v\in\C^n\setminus\{O\}$ such that
$P_\nu(v)=\lambda v$ for some $\lambda\in\C$; the characteristic
direction~$v$ is said {\sl non-degenerate} if $\lambda\ne0$, {\sl
degenerate} otherwise. Then \'Ecalle studied maps~$f$ with at least one
non-degenerate characteristic direction. This is a generic hypothesis,
but not always satisfied: for instance, in the classification we shall
discuss in Section~2 of this note, maps of the classes $(\star^0_1)$,
$(J_1)$ and~$(J_0)$ have {\sl no} non-degenerate characteristic
directions, and so \'Ecalle's methods cannot be applied to them.

We encountered a similar situation in the past studying the
generalization of the Fatou flower theorem to dimension two (see~[A]).
One way of dealing with maps having no non-degenerate characteristic
directions is by blowing up the origin (a technique introduced in this
context by Hakim [H]). In this way we replace $\C^2$ by a manifold~$M$,
the origin by an exceptional divisor~$S$, and the original~$f$ by a
germ~$\tilde f$ of self-map of~$M$ fixing pointwise the exceptional
divisor~$S$. Choosing, as we may, coordinates in which the exceptional
divisor is the line~$\{w_1=0\}$ we can then write
$$
\tilde f(w)=w+w^\nu_1 f^o(w),
\neweq\eqzdue
$$
where $\nu\ge 1$ and ~$f^o$ is a (germ of) holomorphic map not
divisible by~$w_1$. Now, if these coordinates are centered at a
non-degenerate characteristic direction, then necessarily $f^o(O)=O$.
However, there might be other points of the exceptional divisor where
this happens; such points are called {\sl singular points} of~$\tilde f$
(or {\sl singular directions} of~$f$). It turns out that singular
directions always exist, and that they are
the generalization of non-degenerate characteristic directions needed to
get a Fatou flower theorem in dimension~2.

This suggests to study the formal classification of all maps of the
form~\eqzdue\ with $f^o(O)=O$, and this is the aim of this note. Brochero
Martinez has followed a similar approach in~[BM], but he was
looking for normal forms with respect to semi-formal conjugations (that
is conjugations which depends holomorphically on the
$w_2$ variable and formally on the $w_1$ variable), and so his
results are not really comparable to ours. His paper is however
relevant to our study because he also discusses how to deduce information
about the formal classification of the original map~$f$ ([BM,
Theorem~6.2]).\def\autore{Marco Abate and
Francesca Tovena}

In Section~1 of this note we shall describe a general procedure
(inspired by the classical proof of Poincar\'e-Dulac normal forms)
producing normal forms of maps of the form~\eqzdue\ with $f^o(O)=O$ in
dimension~2 (actually, the same procedure works in any dimension; we
restrict ourselves to dimension~2 for simplicity).

In Section~2 we shall apply this procedure to obtain the formal normal
forms of maps satisfying two additional hypotheses. The first one is
just technical: we shall assume that the linear part of~$f^o$ does not
vanish, an hypothesis that can be always satisfied after a finite number
of blow-ups (see~[A]). 
Furthermore, we shall also assume that the linear part of~$f^o$ preserves
the exceptional line~$\{w_1=0\}$. The reason behind this hypothesis is
more conceptual: if $\tilde f$ is the blow-up of a map~$f$, then the
linear part of~$f^o$ (actually, the map~$f^o$ itself) always
preserves the exceptional line~$\{w_1=0\}$ {\it unless} $f$ is
dicritical, that is, writing $f$ as in~\eqzumez, unless $P_\nu(z)=q(z)z$
for a suitable homogeneous polynomial~$q$ of total degree~$\nu-1$, a
very special case. Since the dynamics of maps such that
$f^o$ does not preserve the exceptional line is much easier to study
than the dynamics of maps which do (see~[ABT]), and furthermore Brochero
Martinez [BM, Theorem~6.3] has formal normal forms for maps which are
the blow-up of dicritic maps, we shall restrict ourselvese to map
satisfying this additional hypothesis. So our main
result Theorem~2.1 shall  list the formal normal forms of
2-dimensional maps of the form~\eqzdue\ with
$f^o(O)=O$ and such that the linear part of $f^o$ (is not zero and) preserves the
line~$\{w_1=0\}$.

\smallsect 1. The procedure

We shall denote by $\End(\C^2,O)$ the set of germs at the origin of
holomorphic self-maps of $\C^2$ fixing the origin, and we shall say that
$f\in\End(\C^2,O)$ is {\sl tangent to the identity} if $df_O=\id$.
We are interested in finding formal
normal forms for maps $f\in\End(\C^2,O)$, $f\ne\id_{\C^2}$, tangent to
the identity and satisfying a few additional hypotheses. The first
hypothesis is
\medskip
\item{$({\bf H}_1)$} there exists a curve $S\subset\C^2$ of fixed points
such that $O\in S$ is a regular (i.e., smooth) point of~$S$.
\medskip
\noindent This can always be achieved by blowing-up the origin and
taking as $S$ the exceptional divisor; as discussed in the
introduction, this is an often useful procedure to study the dynamics of
maps tangent to the identity.

Up to a (convergent) change of coordinates, we can assume
that $S=\{z_1=0\}$ near the origin. In particular, all the changes of
coordinates we shall use from now on will preserve the curve
$\{z_1=0\}$.

Since $f$ is tangent to the identity, there exists a maximal
number~$\nu\in\N^*$ (the {\sl order of contact} of $f$ with~$S$; see
[ABT] for an invariant definition in the general case) so that we can
write
$$
f(z)=z+z_1^\nu f^o(z)
$$ 
for a suitable map~$f^o$ not divisible by~$z_1$. Our second hypothesis
then is
\medskip
\item{$({\bf H}_2)$} $f^o(O)=O$, that is the origin is a {\sl
(dynamically) singular} (see [A] and [ABT]) point of~$S$.
\medskip
\noindent The reason for this hypothesis is that, as proved in [A], the
dynamics of $f$ is concentrated around its singular points, and we are
eventually interested in using formal normal forms to study the
dynamics.

We can also expand $f^o$ in series of homogeneous polynomials,
writing
$$
f^o(z)=P^1(z)+P^2(z)+\cdots,
$$
where each $P^j(z)$ is a pair of homogeneous polynomials
in~$z_1$,~$z_2$ of (total) degree~$j$. The least $j\ge 1$ such that
$P^j\not\equiv O$ is (at least in this paper; the definition in~[A] is
slightly different) {\sl pure order} of~$f$ at the origin. 

The third hypothesis we shall (starting from next section) assume is
\medskip
\item{$({\bf H}_3)$} $f$ has pure order~1 at the origin, that is
$P^1\not\equiv O$.
\medskip
\noindent This is a technical hypothesis, used to simplify computations.
However, in [A] it is shown that with a finite number of
blow-ups every $f$ can be transformed in a map with pure order 1 at all
its singular points.

Our last hypothesis will be
\medskip
\item{$({\bf H}_4)$} $P^1$ sends the curve $S=\{z_1=0\}$ into itself. 
\medskip
\noindent As shown in [A] and [ABT], from a dynamical point of view the
most interesting maps tangent to the identity are the ones {\it
tangential} to their fixed point set. We refer to [ABT] for the general
definition of tangential maps; in our setting, $f$ is tangential to $S$
if and only if $f^o$ sends $S$ into itself, that is, writing
$f^o=(f^o_1,f^o_2)$, if $z_1$ divides~$f^o_1$.  So our hypothesis
$({\bf H}_4)$ is a weak version of tangentiality. We shall also discuss
what happens when $f$ is actually tangential to~$S$. 

A few more notations: we shall write $P^j=(P^j_1,P^j_2)$ and
$$
P^j_i(z)=\sum_{k=0}^j a^j_{i,k}z_1^k z_2^{j-k};
$$
in particular, $({\bf H}_4)$ holds if and only if $a^1_{1,0}=0$.
Finally, we shall use the symbol
$O_d$ to indicate any formal series whose expansion in homogeneous polynomials
has no terms of degree less than~$d$. For instance, if $\mu$ is the pure
order of~$f$ at the origin we shall sometimes write
$f^o(z)=P^\mu(z)+O_{\mu+1}$. 

Any formal change of coordinates can be written as the composition of a
linear change of coordinates~$A=\left(\matrix{\alpha_{11}&\alpha_{12}\cr
\alpha_{21}&\alpha_{22}\cr}\right)$ followed by a formal change of
coordinates $\chi$ tangent to the identity, that is whose expansion in
homogeneous polynomials is of the form
$$
\chi(z)=z+H^d(z)+H^{d+1}(z)+\cdots
$$
with $d\ge 2$ and $H^d\not\equiv O$, where again $H^j=(H^j_1,H^j_2)$ is a
pair of homogeneous polynomials of degree~$j$. Furthermore, since
in our setting we are interested only in change of coordinates preserving
the line $\{z_1=0\}$, we must have $\alpha_{12}=0$,
and we should be able to write
$H^j_1=z_1\check H^j_1$ for each $j\ge d$, where $\check H^j_1$ is a
homogeneous polynomial of degree~$j-1$.

Let us now see how a change of coordinates acts on a map satisfying
$({\bf H}_1)$ and $({\bf H}_2)$.

\newthm Proposition \uuno: Let $f\in\End(\C^2,O)$ satisfy hypotheses
$({\bf H}_1)$ and $({\bf H}_2)$, so that it can be written as
$$
f(z)=z+z_1^\nu\bigl(P^\mu(z)+O_{\mu+1}\bigr),
$$
where $\nu\ge 1$ is the order of contact of $f$ with $\{z_1=0\}$, and
$\mu\ge 1$ is the pure order of $f$ at the origin. 
\smallskip
\item{\rm (i)} Let
$$A=\left(\matrix{\alpha_{11}&0\cr
\alpha_{21}&\alpha_{22}\cr}\right)
\neweq\eqlc
$$
be a linear change of coordinates preserving the line $\{z_1=0\}$. Then
$$
A^{-1}\circ f\circ A(z)=z+z_1^\nu\bigl(\alpha_{11}^\nu A^{-1}P^\mu(Az)
	+O_{\mu+1}\bigr).
$$
\item{\rm(ii)}Let
$$
\chi(z)=z+H^d(z)+O_{d+1}
$$
be a formal change of coordinates tangent to the identity and preserving
the line $\{z_1=0\}$, where $d\ge 2$. Then
$$
\eqalign{
\chi^{-1}\circ
f\circ\chi(z)=z+z_1^\nu\bigl[P^\mu(z)&+\cdots+P^{\mu+d-2}(z)\cr
&+P^{\mu+d-1}(z)+\hbox{\rm Jac}(P^\mu)\cdot H^d(z)-\hbox{\rm
Jac}(H^d)\cdot P^\mu(z)+\nu\check H^d_1(z)P^\mu(z)\cr
&{}+O_{\mu+d}\bigr],\cr}
\neweq\eqcoc
$$
\indent where $\hbox{\rm Jac}(H)\cdot v$ denotes the multiplication of
the vector $v\in\C^2$ by the Jacobian matrix of the map~$H$\break\indent
(computed at~$z$).

\pf Part (i) is a trivial computation. To prove \eqcoc, let 
$\chi^{-1}(z)=z+K^e(z)+K^{e+1}(z)+\cdots$ be the expansion of~$\chi^{-1}$
in homogeneous polynomials, with $K^e\not\equiv O$. From
$\chi\circ\chi^{-1}=\id$ we get
$$
\sum_{j\ge e} K^j(z)+H^d\left(z+\sum_{j\ge e}
K^j(z)\right)+O_{d+1}(z)\equiv O.
$$
Replacing $z$ by $\lambda z$, with $\lambda\in\C^*$, we get
$$
\lambda^e \sum_{j\ge e}\lambda^{j-e}
K^j(z)+\lambda^d H^d\left(z+\sum_{j\ge e}\lambda^{j-1}
K^j(z)\right)+\lambda^{d+1} g(z)\equiv O
$$
for a suitable (formal) map $g$. Dividing by $\lambda^d$ and then
letting $\lambda\to 0$ we then find $e=d$ and $K^e=-H^d$. In particular,
then,
$$
\hbox{\rm Jac}(\chi^{-1})=I-\hbox{\rm Jac}(H^d)+O_d,
\neweq\eqchimeno
$$
where $I$ is the identity matrix, and any partial derivative
of~$\chi^{-1}_i$ of order~$2\le\ell\le d$ is an~$O_{d-\ell}$.

Now,
$$
\eqalign{
\chi^{-1}\circ f\circ\chi(z)&=\chi^{-1}\Bigl(\chi(z)+\chi_1(z)^\nu
f^o\bigl(\chi(z)\bigr)\Bigr)\cr
&=z+z_1^\nu\check\chi_1(z)^\nu\hbox{\rm Jac}(\chi^{-1})\cdot
f^o\bigl(\chi(z)\bigr)\cr &\qquad
+{1\over2}z_1^{2\nu}\check\chi_1(z)^{2\nu}\sum_{h,k=1}^2
{\de^2\chi^{-1}\over\de z_h\de z_k}\bigl(\chi(z)\bigr) 
f^o_h\bigl(\chi(z)\bigr)f^o_k\bigl(\chi(z)\bigr)+\cdots\cr
&=z+z_1^\nu\left[\check\chi_1(z)^\nu \hbox{\rm Jac}(\chi^{-1})\cdot
f^o\bigl(\chi(z)\bigr)+O_{\mu+d}\right].
\cr}
\neweq\eqconto
$$
Putting
$$
\check\chi_1(z)^\nu=\bigl(1+\check H^d_1(z)+O_d\bigr)^\nu=1+\nu
\check H^d_1(z)+O_d,
$$
$$
P^j\bigl(\chi(z)\bigr)=P^j\bigl(z+H^d(z)+O_{d+1}\bigr)=P^j(z)+
\hbox{\rm Jac}(P^j)\cdot H^d(z)+O_{d+j},
$$
and \eqchimeno\ in \eqconto\ we finally get \eqcoc.\qedn

Formula \eqcoc\ suggests the introduction of a family of linear maps.
Let $\C_d[z_1,z_2]$ denote the space of homogeneous polynomials of
degree~$d$ in two variables, and set $\tilde{\ca
V}_d=\C_d[z_1,z_2]^2$ (that is, $\tilde{\ca V}_d$ is the space of pairs
of homogeneous polynomials of degree~$d$ in two variables) and
$$
\ca V_d=\{H\in\tilde{\ca V}_d\mid H_1(0,\cdot)\equiv 0\},
$$
so that $H\in\ca V_d$ implies that there is $\check
H_1\in\C_{d-1}[z_1,z_2]$ such that $H_1\equiv z_1\check H_1$.

Given $P\in\tilde{\ca V}_\mu$, for any $d\ge 2$ we can then define a
linear map $L_{P,d}\colon\ca V_d\to\tilde{\ca V}_{d+\mu-1}$
by setting
$$
L_{P,d}(H)=\hbox{\rm Jac}(P)\cdot H-\hbox{\rm
Jac}(H)\cdot P+\nu\check H_1P.
$$
We explicitely remark that if $P\in\ca V_\mu$ then $L_{P,d}$ sends $\ca
V_d$ into~$\ca V_{d+\mu-1}$. 

From formula \eqcoc\ we then deduce the following:

\newthm Corollary \udue: Let $f\in\End(\C^2,O)$ of the form
$$
f(z)=z+z_1^\nu\bigl(P^\mu(z)+\cdots+P^{\mu+d-1}(z)+O_{\mu+d}\bigr),
$$
and take $Q\in\tilde{\ca V}_{\mu+d-1}$. 
Then we can find a change of coordinates of the form $\chi(z)=z+H^d(z)$
with $H^d\in\ca V_d$ such that
$$
\chi^{-1}\circ
f\circ\chi(z)=z+z_1^\nu\bigl(P^\mu(z)+\cdots+P^{\mu+d-1}(z)+Q(z)+O_{d+\mu}\bigr)
$$
if and only if $Q$ belongs to the image of $L_{P^\mu,d}$. In particular,
we can find a change of coordinates of the form~$\chi(z)=z+H^d(z)$
with $H^d\in\ca V_d$ such that
$$
\chi^{-1}\circ
f\circ\chi(z)=z+z_1^\nu\bigl(P^\mu(z)+\cdots+P^{\mu+d-2}(z)+O_{d+\mu}\bigr)
$$
if and only if $P^{\mu+d-1}$ belongs to the image of $L_{P^\mu,d}$.

This corollary suggests the announced procedure for finding formal normal forms. First
of all, one uses a linear change of variables to put~$P^\mu$ in normal
form. Then one uses a change of variables of the form $\chi(z)=z+H^2(z)$
to put $P^{\mu+1}$ in normal form, subtracting elements of the
image of~$L_{P^\mu,2}$. Then one proceeds by induction: with a change
of variables of the form $\chi(z)=z+H^d(z)$, subtracting
elements of the image of~$L_{P^\mu,d}$ one can put $P^{\mu+d-1}$
in normal form without modifying the preceding terms. Composing all these changes of
variables for $d\to\infty$ we get a formal change of variables putting
the map in formal normal form.

It is worthwhile to remark that with this procedure it is possible to
obtain normal forms up to any given finite order by using
{\it polynomials} change of coordinates. This can be useful, because
in the study of the dynamics in a neighbourhood of the origin higher
order terms are often negligible. 

In the next section we shall apply this procedure to
maps satisfying hypotheses $({\bf H}_3)$ and $({\bf H}_4)$. 

\smallsect 2. Normal forms

From now on we shall assume hypotheses $({\bf H}_1)$, $({\bf H}_2)$ and
$({\bf H}_3)$; so let $f\in\End(\C^2,O)$ be of the form
$$
f(z)=z+z_1^\nu\bigl(P^1(z)+P^2(z)+\cdots)
$$
with $P^1\not\equiv O$. The first step in our procedure consists in
putting~$P^1$ in normal form using a linear change of coordinates~$A$ of
the form~\eqlc. Using the notations introduced in the previous section,
Proposition~\uuno.(i) says that $P^1$ is sent by~$A$ in
$$
\alpha_{11}^\nu A^{-1}\circ P^1\circ A=\alpha_{11}^\nu\left|
\matrix{a^1_{1,1}+{\alpha_{21}\over\alpha_{11}}a^1_{1,0}&
{\alpha_{22}\over\alpha_{11}}a^1_{1,0}\cr
{\alpha_{21}\over\alpha_{22}}(a^1_{2,0}-a^1_{1,1})+{\alpha_{11}\over
\alpha_{22}}a^1_{2,1}-{\alpha^2_{21}\over\alpha_{11}\alpha_{22}}
a^1_{1,0}&a^1_{2,0}-{\alpha_{21}\over\alpha_{11}}a^1_{1,0}
\cr}
\right|.
$$
We have several cases to consider.
\smallskip
\item{--} If $a^1_{1,0}\ne 0$ (and thus, in particular, $({\bf H}_4)$
does not hold), we can choose $A$ so that
$$
\alpha^{\nu-1}_{11}\alpha_{22} a^1_{1,0}=1\qquad\hbox{and}\qquad
a^1_{2,0}-{\alpha_{21}\over\alpha_{11}}a^1_{1,0}=0.
$$
Then
$$
\alpha^\nu_{11}\left(a^1_{1,1}+{\alpha_{21}\over\alpha_{11}}a^1_{1,0}
\right)=\alpha^\nu_{11}\hbox{tr}(P^1)
$$
and
$$
\alpha^\nu_{11}\left({\alpha_{21}\over\alpha_{22}}(a^1_{2,0}-a^1_{1,1})
+{\alpha_{11}\over\alpha_{22}}a^1_{2,1}-{\alpha^2_{21}\over\alpha_{11}
\alpha_{22}}a^1_{1,0}\right)=-a^{2\nu}_{11}\det(P^1).
$$
So we have three possibilities:
\smallskip
\itemitem{$(T_\lambda)$} If $\hbox{tr}(P^1)\ne 0$ we can choose $A$ so
that $\alpha^\nu_{11}\hbox{tr}(P^1)=1$ and we have reduced~$P^1$ to
$$
\left|\matrix{1&1\cr \lambda&0\cr}\right|,
$$
where $\lambda=-\det(P^1)/\hbox{tr}(P^1)^2\in\C$.
\smallskip
\itemitem{$(N_1)$} If instead $\hbox{tr}(P^1)=0$ but $\det(P^1)\ne 0$, we
can choose $A$ so that $-a^{2\nu}_{11}\det(P^1)=1$, and we have
reduced~$P^1$ to
$$
\left|\matrix{0&1\cr 1&0\cr}\right|.
$$
\smallskip
\itemitem{$(N_0)$} Finally, if $\hbox{tr}(P^1)=0=\det(P^1)$, we
have reduced~$P^1$ to
$$
\left|\matrix{0&1\cr 0&0\cr}\right|.
$$
\medskip
\item{--} On the other hand, if $a^1_{1,0}=0$ (and so $({\bf
H}_4)$ holds) and $a^1_{2,0}\ne 0$, we can choose $A$ so that~$a_{11}^\nu
a^1_{2,0}=1$, and we again have three subcases to consider:
\smallskip
\itemitem{$(\star_1^\lambda)$} If $a^1_{1,1}\ne a^1_{2,0}$ we can
choose~$A$ so that ${\alpha_{21}\over\alpha_{22}}(a^1_{2,0}-a^1_{1,1})+{\alpha_{11}\over
\alpha_{22}}a^1_{2,1}=0$, and we have reduced $P^1$ to
$$
\left|\matrix{\lambda&0\cr 0&1\cr}\right|,
$$
where $\lambda=a^1_{1,1}/a^1_{2,0}\ne 1$ is the {\sl residual index} of
$f$ at the origin along~$\{z_1=0\}$ (see [A] for the definition of
residual index).
\smallskip
\itemitem{$(J_1)$} If instead $a^1_{1,1}=a^1_{2,0}$ but $a^1_{2,1}\ne
0$ we can choose~$A$ so that
${\alpha_{11}^{\nu+1}\over
\alpha_{22}}a^1_{2,1}=1$, and we have reduced $P^1$ to
$$
\left|\matrix{1&0\cr 1&1\cr}\right|.
$$
\smallskip
\itemitem{$(\star_1^1)$} Finally, if $a^1_{1,1}=a^1_{2,0}$ and
$a^1_{2,1}=0$ we have reduced $P^1$ to
$$
\left|\matrix{1&0\cr 0&1\cr}\right|.
$$
\medskip
\item{--} And now, if $a^1_{1,0}=a^1_{2,0}=0$ we have two last cases to
consider:
\smallskip
\itemitem{$(\star_2)$} If $a_{1,1}^1\ne 0$ we can choose $A$ so that
$\alpha^\nu_{11}a^1_{1,1}=1$ and 
${\alpha_{11}\over
\alpha_{22}}a^1_{2,1}-{\alpha_{21}\over\alpha_{22}}a^1_{1,1}=0$, and we
have reduced~$P^1$ to
$$
\left|\matrix{1&0\cr 0&0\cr}\right|.
$$
\smallskip
\itemitem{$(J_0)$} Finally, if $a_{1,1}^1=0$ also, we can choose $A$ so
that
${\alpha_{11}^{\nu+1}\over\alpha_{22}}a^1_{2,1}=1$, and we
have reduced~$P^1$ to
$$
\left|\matrix{0&0\cr 1&0\cr}\right|.
$$
\medskip
\noindent So we have eight possible normal forms for the linear part of
a map satisfying $({\bf H}_1)$, $({\bf H}_2)$ and $({\bf H}_3)$, and five
for the linear part of a map satisfying also $({\bf H}_4)$.

To apply the procedure described at the end of the previous
section, we must study the linear maps $L_{P^\mu,d}\colon\ca
V_d\to\tilde{\ca V}_{d+\mu-1}$. For the moment, we shall not use
hypotheses~$({\bf H}_3)$ and $({\bf H}_4)$, specializing to~$\mu=1$
and $a^1_{1,0}=0$ only later.

We shall use monomials as a basis of $\tilde{\ca V}_d$: setting
$$
v_d^h=\cases{(z_1^h z_2^{d-h},0)&for $h=0,\ldots,d$,\cr
\noalign{\smallskip}
(0,z_1^k z_2^{d-k})&for $h=d+1,\ldots,2d+1$, where $k=h-(d+1)$,\cr}
$$
we see that $\{v_d^0,\ldots,v_d^{2d+1}\}$ is a basis of $\tilde{\ca
V}_d$, while $\{v_d^1,\ldots,v_d^{2d+1}\}$ is a basis of~$\ca V_d$. 

Using the notations introduced in the previous section, we can write
$$
P^\mu(z)=\left(\sum_{i=0}^\mu a^\mu_{1,i}z_1^iz_2^{\mu-i},
\sum_{i=0}^\mu a^\mu_{2,i}z_1^iz_2^{\mu-i}\right)=
\sum_{i=0}^\mu a^\mu_{1,i} v_\mu^i+\sum_{i=0}^\mu a^\mu_{2,i}
v_\mu^{d+1+i}.
$$
Then it is not difficult to compute the action of $L_{P^\mu,d}$ on the
elements of the basis of~$\ca V_d$:
$$
L_{P^\mu,d}(v_d^h)=\sum_{i=0}^{\mu+1}
[(\nu+i-h)a^\mu_{1,i}-(d-h)a^\mu_{2,i-1}]\,v^{h+i-1}_{d+\mu-1}+
\sum_{i=0}^\mu (\nu+i)a^\mu_{2,i}\, v_{d+\mu-1}^{h+i+d}
$$
for $h=1,\ldots,d$, where we have put $a^\mu_{1,\mu+1}=0=a^\mu_{2,-1}$, 
and
$$
L_{P^\mu,d}(v_d^{k+d+1})=\sum_{i=0}^{\mu-1}
(\mu-i)a^\mu_{1,i}\,v^{k+i}_{d+\mu-1}+
\sum_{i=0}^{\mu+1}[(\mu+k+1-d-i)a^\mu_{2,i-1}-k a^\mu_{1,i}]\,
v_{d+\mu-1}^{k+i+d}
$$
for $k=0,\ldots,d$, where again we have put
$a^\mu_{1,\mu+1}=0=a^\mu_{2,-1}$.

Specializing to $\mu=1$ we get
$$
L_{P^1,d}(v_d^h)=\sum_{i=0}^{2}
[(\nu+i-h)a^1_{1,i}-(d-h)a^1_{2,i-1}]\,v^{h+i-1}_d+
\sum_{i=0}^1 (\nu+i)a^1_{2,i}\, v_d^{h+i+d}
\neweq\eqLuno
$$
for $h=1,\ldots,d$,  
and
$$
L_{P^1,d}(v_d^{k+d+1})=
a^1_{1,0}\,v^k_d+
\sum_{i=0}^{2}[(k+2-d-i)a^1_{2,i-1}-k a^1_{1,i}]\,
v_d^{k+i+d}
\neweq\eqLdue
$$
for $k=0,\ldots,d$.

When $({\bf H}_3)$ and $({\bf H}_4)$ hold we have $\mu=1$ and
$a^1_{1,0}=0$ in \eqLuno\ and~\eqLdue; in particular, the
image of $L_{P^1,d}$ is always contained in~$\ca V_d$. This implies that
we cannot ever remove from the normal form monomials proportionals to
$(z_2^d,0)$; however, these monomials appear if and only if
$f$ is not tangential.

We shall assume now that $P^1$ is in normal form, and apply our
procedure.
\medskip
$(\star_1^\lambda)$. In this case we have $a^1_{1,1}=\lambda$,
$a^1_{2,0}=1$ and $a^1_{1,0}=a^1_{2,1}=0$, and thus \eqLuno\ and
\eqLdue\ become
$$
L_{P^1,d}(v_d^h)=
[(\nu+1-h)\lambda-(d-h)]\,v^h_d+\nu\, v_d^{h+d}
$$
for $h=1,\ldots,d$,  
and
$$
L_{P^1,d}(v_d^{k+d+1})=
[(k+1-d)-k\lambda]\,v_d^{k+1+d}
$$
for $k=0,\ldots,d$. The matrix representing $L_{P^1,d}$ with respect to
the basis
$\{v^1_d,\ldots,v^{2d+1}_d\}$ of~$\ca V_d$ is lower triangular; in
particular, $L_{P^1,d}$ is surjective if and only if all the elements on
the diagonal are different from zero. Now, these elements are of two
forms: 
$$
\sigma_{d,h}=(\nu+1-h)\lambda-(d-h)
$$
with $h=1,\ldots,d$, and
$$
\tau_{d,k}=(k+1-d)-k\lambda
$$
with $k=0,\ldots,d$. In particular, $\sigma_{\nu+1,\nu+1}=0$ always;
this implies that we cannot ever remove from the normal form the
monomial $(z_1^{\nu+1},0)$.  

Let $E\subset\C$ denote the set of $\lambda\in\C$ such that 
$\sigma_{d,h}=0$ for some pair $(d,h)\ne(\nu+1,\nu+1)$, or such that
$\tau_{d,k}=0$ for some pair $(d,k)$. Clearly, if $\lambda\notin E$ we
have that $L_{P^1,d}$ is surjective for all $d\ne\nu+1$, and
$\Im(L_{P^1,\nu+1})=\hbox{Span}(v^1_{\nu+1},\ldots,v^\nu_{\nu+1},
v^{\nu+2}_{\nu+1},\ldots,v^{2\nu+3}_{\nu+1})$, that is we can remove
from the normal form {\it all} the monomials {\it but} $(z_1^{\nu+1},0)$
and $(z_2^d,0)$.
This means that if $\lambda\notin E$ the formal normal form of $f$ is
$$
\hat f(z)=\bigl(z_1+z_1^\nu[\lambda z_1+a z_1^{\nu+1}+z_2^2t(z_2)],
z_2+z_1^\nu z_2\bigr)
$$
with $a\in\C$ and $t\in\C[\![\zeta]\!]$, where $\C[\![\zeta]\!]$ is
the space of formal power series in 1 variable, with
$t\equiv 0$ if and only if
$f$ is tangential.

So it becomes important to determine~$E$. First of all, $\tau_{d,k}=0$
if and only if $\lambda=-(d-k-1)/k$ with $1\le k\le d$ and $d\ge 2$, and
it is easy to see that
$$
\left\{-{d-k-1\over k}\biggm| 1\le k\le d, d\ge 2\right\}=\left\{
{1\over d}\biggm| d\ge 2\right\}\cup\{0\}\cup\Q^-.
$$
On the other hand, $\sigma_{d,h}=0$ if and only if
$\lambda=(d-h)/(\nu+1-h)$ with $1\le h\le d$, $h\ne\nu+1$, $d\ge 2$. Now,
if $h>\nu+1$ then we get a negative rational number, already considered
before. For the same reason we can for the moment disregard the case
$h=d$, and limit ourselves to the case $1\le h<\min\{\nu+1,d\}$, when
$\lambda$ is a positive rational number. Setting $q=\nu+1-h$ and
$p=d-h=d+q-(\nu+1)$, we see that we must have
$\max\{0,\nu+1-d\}<q\le\nu$ and $\max\{d-(\nu+1),0\}<p\le d-1$.
Since $d\ge 2$ is generic, we get
$$
\left\{{d-h\over\nu+1-h}\biggm|1\le h\le d, h\ne\nu+1, d\ge
2\right\}\cap\Q^+=\bigcup_{q=1}^\nu {1\over q}\N,
$$
and so
$$
E=\bigcup_{q=1}^\nu {1\over q}\N\cup\left\{
{1\over d}\biggm| d\ge 2\right\}\cup\{0\}\cup\Q^-.
$$
We now study $L_{P^1,\mu}$ when $\lambda\in E$.
\medskip
\item{--} $\lambda=0$. In this case $\tau_{d,k}=0$ if and only if
$k=d-1$, and $\sigma_{d,h}=0$ if and only if $h=d$. So $v^{k+d+1}_d\in
\Im L_{P^1,d}$ for $k=0,\ldots,d-2,d$, and thus $v^h_d\in
\Im L_{P^1,d}$ for $h=1,\ldots,d-1$. Furthermore, $\nu v^{2d}_d=L_{P^1,d}
(v^d_d)$, and so $\Im
L_{P^1,d}=\hbox{Span}(v^1_d,\ldots,v^{d-1}_d,v^{d+1}_d,\ldots,
v^{2d+1}_d)$. Therefore the formal normal form of $f$ in this case is
$$
\hat f(z)=\bigl(z_1+z_1^\nu[z_1^2 g(z_1)+z_2^2 t(z_2)], z_2+z_1^\nu z_2
\bigr),
$$
where $g$,~$t\in\C[\![\zeta]\!]$, with $t\equiv 0$ if and only if $f$ is
tangential.
\smallskip
\item{--}$\lambda=-p/q\in\Q^-$, with $(p,q)=1$. In this case we have
$\tau_{d,k}=0$ if and only if $d-1-k=kp/q$; in particular, since $p$ and
$q$ are relatively prime, we must have $k=ql$ and $d-1-k=pl$ for some
$l\ge 1$. 

Now, if $\tau_{d,k}\ne 0$ then $v^{k+d+1}_d\in\Im L_{P^1,d}$. If
$\tau_{d,k}=0$ with $1\le k\le d-2$, we have
$$
L_{P^1,d}(v^{k+1}_d)=-\lambda\nu v^{k+1}_d+\nu v^{k+d+1}_d;
$$
\item{}hence ($\sigma_{d,k+1}\ne0$ and) we
can remove from the normal form monomials of the form
$v^{k+1+d}_d$ if we retain monomials of the form
$v^{k+1}_d=(z_1^{k+1}z_2^{d-k-1},0)$. Notice that we can write
$$
z_1^{k+1}z_2^{d-k-1}=z_1(z_1^q z_2^p)^l.
$$
\item{}On the other hand, $\sigma_{d,h}=0$ if and only if
$d-h=(h-\nu-1)p/q$; in particular we must have $h-\nu-1=lq$ and $d-h=lp$
for some $l\ge 1$. This means that $v^h_d\in\Im L_{P^1,d}$ except when we
can write
$$
z_1^hz_2^{d-h}=z_1^{\nu+1}(z_1^qz_2^p)^l.
$$
Summing up, we have shown that the formal normal form in this case is
$$
\hat f(z)=\bigl(z_1+z_1^\nu[z_1g_1(z_1^qz_2^p)+z_1^{\nu+1}g_2(z_1^qz_2^p)
+z_2^2t(z_2)],z_2+z_1^\nu z_2\bigr)
$$
where $g_1$, $g_2$, $t\in\C[\![\zeta]\!]$ with $g_1(0)=\lambda$, and with
$t\equiv 0$ if and only if $f$ is tangential.
\smallskip
\item{--}$\lambda=1$. In this case we have $\tau_{d,k}\ne 0$ always, and
$\sigma_{d,h}=0$ if and only if $d=\nu+1$. This means that $\Im
L_{P^1,d}=\ca V_d$ if $d\ne\nu+1$, and $\Im
L_{P^1,\nu+1}=\hbox{Span}(v^{\nu+2}_{\nu+1},\ldots,v^{2\nu+3}_{\nu+1})$,
and so the normal form in this case is
$$
\hat f(z)=\bigl(z_1+z_1^\nu[z_1+z_1 p_\nu(z_1,z_2)+z_2^2t(z_2)],z_2+
z_1^\nu z_2\bigr)
$$
where $p_\nu\in\C_\nu[z_1,z_2]$ and $t\in\C[\![\zeta]\!]$, with
$t\equiv 0$ if and only if $f$ is tangential.
\smallskip
\item{--}$\lambda=1/q$, with $q\ge 2$. In this case $\tau_{d,k}=0$ if
and only if $d=k=q$. On the other hand, $\sigma_{d,h}=0$ if and only if
$\nu+1-h=(d-h)q$; in particular, if we set $l=d-h$ we must have
$h=\nu+1-lq$ with $1\le l<(\nu+1)/q$, and thus we can write
$$
z_1^hz_2^{d-h}=z_1^{\nu+1}(z_1^{-q}z_2)^l.
$$ 
Arguing as before we then see that the formal normal form in this case is
$$
\hat f(z)=\bigl(z_1+z_1^\nu[\lambda
z_1+z_1^{\nu+1}p_o(z_1^{-q}z_2)+z_2^2t(z_2)],
z_2+z_1^\nu(z_2+az_1^q)\bigr),
$$
where $a\in\C$, $p_o\in\C[\zeta]$ is a polynomial of degree at most
$\nu/q$ (and in particular it is just a constant if~$q\ge\nu+1$),
and $t\in\C[\![\zeta]\!]$, with $t\equiv 0$ if and only if $f$ is
tangential.
\smallskip
\item{--}$\lambda=p/q$, with $1\le q\le\nu$, $p\ge 2$ and $(p,q)=1$. In
this case $\tau_{d,k}\ne 0$ always, and $\sigma_{d,h}=0$ if and only if
$d-h=(\nu+1-h)p/q$. Therefore we can write $\nu+1-h=lq$, $d-h=lp$ and
$$
z_1^hz_2^{d-h}=z_1^{\nu+1}(z_1^{-q}z_2^p)^l,
$$
with $1\le l\le\nu/q$. Hence arguing as before we see that the formal
normal form in this case is
$$
\hat f(z)=\bigl(z_1+z_1^\nu[\lambda z_1+z_1^{\nu+1}p_o(z_1^{-q}z_2^p)
+z_2^2t(z_2)], z_2+z_1^\nu z_2\bigr),
$$
where $p_o\in\C[\zeta]$ is a polynomial of degree at most
$\nu/q$ (and in particular it is just a constant if~$q\ge\nu+1$),
and $t\in\C[\![\zeta]\!]$, with $t\equiv 0$ if and only if $f$ is
tangential.
\medskip
\noindent This ends the discussion of the case $(\star_1^\lambda)$. We
can now deal with the remaining cases:
\medskip
$(J_1)$. In this case we have $a^1_{1,1}=a^1_{2,0}=a^1_{2,1}=1$ and
$a^1_{1,0}=0$, and thus \eqLuno\ and \eqLdue\ become
$$
L_{P^1,d}(v_d^h)=(\nu+1-d)v^h_d-(d-h)v_d^{h+1}
+\nu v^{h+d}_d+(\nu+1) v_d^{h+d+1}
$$
for $h=1,\ldots,d$,  
and
$$
L_{P^1,d}(v_d^{k+d+1})=(1-d)v_d^{k+d+1}-(d-k)v_d^{k+d+2}
$$
for $k=0,\ldots,d$. The matrix representing $L_{P^1,d}$ with respect to
the basis
$\{v^1_d,\ldots,v^{2d+1}_d\}$ of~$\ca V_d$ is lower triangular; in
particular, $L_{P^1,d}$ is surjective if and only if $d\ne\nu+1$.
Furthermore, it is not difficult to check that $\Im
L_{P^1,\nu+1}=\hbox{Span}(v^2_{\mu+1},\ldots,v^{2\mu+3}_{\mu+1})$, and
thus the formal normal form in this case is
$$
\hat f(z)=\bigl(z_1+z_1^\nu[z_1+az_1z_2^\nu+z_2^2 t(z_2)], z_2+z_1^\nu
(z_1+z_2)\bigr),
$$
where $a\in\C$ and $t\in\C[\![\zeta]\!]$, with $t\equiv 0$ if and only if
$f$ is tangential.

\medskip
$(\star_2)$. In this case we have $a^1_{1,1}=1$ and
$a^1_{2,0}=a^1_{2,1}=a^1_{1,0}=0$, and thus \eqLuno\ and \eqLdue\ become
$$
L_{P^1,d}(v_d^h)=(\nu+1-h)\,v^h_d
$$
for $h=1,\ldots,d$,  
and
$$
L_{P^1,d}(v_d^{k+d+1})=
-k v_d^{k+d+1}
$$
for $k=0,\ldots,d$. The matrix representing $L_{P^1,d}$ with respect to
the basis
$\{v^1_d,\ldots,v^{2d+1}_d\}$ of~$\ca V_d$ is diagonal; in
particular, $L_{P^1,d}(v^h_d)=O$ if and only if $h=\nu+1$ (and
$d\ge\nu+1$), and $L_{P^1,d}(v^{k+d+1}_d)=O$ if and only if $k=0$.
Therefore the formal normal form in this case is
$$
\hat f(z)=\bigl(z_1+z_1^\nu[z_1+z_1^{\nu+1}g_1(z_2)+z_2^2t(z_2)],
z_2+z_1^\nu z_2^2g_2(z_2)\bigr),
$$
where $g_1$, $g_2$, $t\in\C[\![\zeta]\!]$, where $t\equiv 0$ if and only
if $f$ is tangential.

\medskip
$(J_0)$. In this case we have $a^1_{2,1}=1$ and
$a^1_{2,0}=a^1_{1,1}=a^1_{1,0}=0$, and thus \eqLuno\ and \eqLdue\ become
$$
L_{P^1,d}(v_d^h)=-(d-h)\,v^{h+1}_d+
(\nu+1)\, v_d^{h+d+1}
$$
for $h=1,\ldots,d$,  
and
$$
L_{P^1,d}(v_d^{k+d+1})=
-(d-k)\,v_d^{k+d+2}
$$
for $k=0,\ldots,d$. The matrix representing $L_{P^1,d}$ with respect to
the basis
$\{v^1_d,\ldots,v^{2d+1}_d\}$ of~$\ca V_d$ is strictly lower triangular,
but with no zeroes on the diagonal just below the main diagonal. It
follows easily that the image of~$L_{P^1,d}$ is generated by
$\{v^2_d,\ldots,v^d_d,v^{d+2}_d,\ldots,v^{2d+1}_d\}$, and thus the
formal normal form in this case is
$$
\hat f(z)=\bigl(z_1+z_1^\nu[z_1z_2 g_1(z_2)+z_2^2t(z_2)],
z_2+z_1^\nu[z_1+z_2^2 g_2(z_2)]\bigr),
$$
where $g_1$, $g_2$, $t\in\C[\![\zeta]\!]$, where $t\equiv 0$ if and only
if $f$ is tangential.
\medbreak

Summing up, we have proved the following

\newthm Theorem \duno: Let $f\in\End(\C^2,O)$, $f\ne\id_{\C^2}$, be a
map tangent to the identity fixing pointwise the line~$S=\{z_1=0\}$.
Assume moreover that the origin is a singular point for $f$ on~$S$, 
that $f$ has order of contact $\nu\ge 1$ with~$S$, pure order~$1$ at the
origin, and that
$({\bf H}_4)$ holds. Let
$$
E=\bigcup_{q=1}^\nu {1\over q}\N\cup\left\{
{1\over d}\biggm| d\ge 2\right\}\cup\{0\}\cup\Q^-.
$$
Then $f$ is formally conjugated to one (and only
one) of the following normal forms:
\smallskip
\item{$(\star_1^\lambda)$} with $\lambda\in\C\setminus E$: 
$$
\hat f(z)=\bigl(z_1+z_1^\nu[\lambda z_1+a z_1^{\nu+1}+z_2^2t(z_2)],
z_2+z_1^\nu z_2\bigr)
$$
where $a\in\C$ and $t\in\C[\![\zeta]\!]$ is a generic power
series;\item{$(\star_1^\lambda)$} with $\lambda=0$: 
$$
\hat f(z)=\bigl(z_1+z_1^\nu[z_1^2 g(z_1)+z_2^2 t(z_2)], z_2+z_1^\nu z_2
\bigr),
$$
where $g$, $t\in\C[\![\zeta]\!]$ are generic power series;
\item{$(\star_1^\lambda)$} with $\lambda=-p/q\in\Q^-$ and $(p,q)=1$:
$$
\hat f(z)=\bigl(z_1+z_1^\nu[z_1g_1(z_1^qz_2^p)+z_1^{\nu+1}g_2(z_1^qz_2^p)
+z_2^2t(z_2)],z_2+z_1^\nu z_2\bigr)
$$
where $g_1$,~$g_2$, $t\in\C[\![\zeta]\!]$ are generic power series with
$g_1(0)=\lambda$;
\item{$(\star_1^\lambda)$} with $\lambda=1/q$ and $q\in\N$, $q\ge
2$: 
$$
\hat f(z)=\bigl(z_1+z_1^\nu[\lambda
z_1+z_1^{\nu+1}p_o(z_1^{-q}z_2)+z_2^2t(z_2)],
z_2+z_1^\nu(z_2+az_1^q)\bigr),
$$
where $a\in\C$, $p_o\in\C[\zeta]$ is a polynomial of degree at most
$\nu/q$, and $t\in\C[\![\zeta]\!]$ is a generic power
series;
\item{$(\star_1^\lambda)$} with $\lambda=p/q\in\Q$, $1\le
q\le\nu$, $p\ge 2$ and $(p,q)=1$: 
$$
\hat f(z)=\bigl(z_1+z_1^\nu[\lambda z_1+z_1^{\nu+1}p_o(z_1^{-q}z_2^p)
+z_2^2t(z_2)], z_2+z_1^\nu z_2\bigr),
$$
where $p_o\in\C[\zeta]$ is a polynomial of degree at most~$\nu/q$, and  
$t\in\C[\![\zeta]\!]$ is a generic power series;
\item{$(\star_1^\lambda)$} with $\lambda=1$: 
$$
\hat f(z)=\bigl(z_1+z_1^\nu[z_1+z_1 p_\nu(z_1,z_2)+z_2^2t(z_2)],z_2+
z_1^\nu z_2\bigr)
$$
where
$p_\nu\in\C_\nu[z_1,z_2]$ is a generic homogeneous polynomial of
degree~$\nu$, and  $t\in\C[\![\zeta]\!]$ is a generic power
series;
\item{$(J_1)$}
$$
\hat f(z)=\bigl(z_1+z_1^\nu[z_1+az_1z_2^\nu+z_2^2 t(z_2)], z_2+z_1^\nu
(z_1+z_2)\bigr),$$
where $a\in\C$ and $t\in\C[\![\zeta]\!]$ is a generic power
series;
\item{$(\star_2)$} 
$$
\hat f(z)=\bigl(z_1+z_1^\nu[z_1+z_1^{\nu+1}g_1(z_2)+z_2^2 t(z_2)],
z_2+z_1^\nu z_2^2 g_2(z_2)\bigr),
$$
where $g_1$,~$g_2$, $t\in\C[\![\zeta]\!]$ are generic power series;
\item{$(J_0)$}
$$
\hat f(z)=\bigl(z_1+z_1^\nu[z_1z_2
g_1(z_2)+z_2^2t(z_2)],z_2+z_1^\nu[z_1+z_2^2 g_2(z_2)]
\bigr),
$$
 where $g_1$,~$g_2$, $t\in\C[\![\zeta]\!]$ are generic power
series.
\smallskip
\noindent In all these cases, $t\equiv 0$ if and only if $f$ is
tangential.

\bigbreak
\noindent {\bf References.}\setref{ABT}
\medskip
\art A M. Abate: The residual index and the dynamics of holomorphic
maps tangent to the identity! Duke Math. J.! 107 2001 173-207

\art ABT M. Abate, F. Bracci, F. Tovena: Index theorems for
holomorphic self-maps!   Ann. of Math.! 159 2004 819-864

\art BM F.E. Brochero Mart\'\i nez: Groups of germs of
analytic diffeomorphisms in $(\C^2,O)$! J. Dynamic. Control Systems! 9
2003 1-32

\art C C. Camacho: On the local structure of conformal mappings
and holomorphic vector fields! Ast\'e\-risque! 59--60 1978 83-94

\book \'E1 J. \'Ecalle: Les fonctions r\'esurgentes. Tome
I: Les alg\`ebres de fonctions r\'esurgentes! Publ. Math. Orsay {\bf
81-05,} Universit\'e de Paris-Sud, Orsay, 1981 

\book \'E2 J. \'Ecalle: Les fonctions r\'esurgentes. Tome
I\negthinspace I: Les fonctions r\'esurgentes appliqu\'ees
\`a l'it\'eration! Publ. Math. Orsay {\bf 81-06,} Universit\'e de
Paris-Sud, Orsay, 1981 

\book \'E3 J. \'Ecalle: Les fonctions r\'esurgentes. Tome I\negthinspace
I\negthinspace I: L'\'equation du pont et la classification analytique
des objects locaux! Publ. Math. Orsay {\bf 85-05,} Universit\'e de
Paris-Sud, Orsay, 1985 

\coll \'E4 J. \'Ecalle: Iteration and analytic classification of local
diffeomorphisms of $\C^\nu$! Iteration theory and its
functional equations! Lect. Notes in Math. {\bf 1163,} Springer-Verlag,
Berlin, 1985, pp. 41--48

\art H M. Hakim: Analytic transformations of $(\C^p,0)$ tangent to the
identity! Duke Math. J.! 92 1998 403-428

\coll I Yu.S. Il'yashenko: Nonlinear Stokes phenomena! Nonlinear Stokes
phenomena! Adv. in Soviet Math. {\bf 14,} Am. Math. Soc., Providence,
1993, pp.~1--55

\art M1 B. Malgrange: Travaux d'\'Ecalle et de Martinet-Ramis sur les
syst\`emes  dynamiques! Ast\'erisque! 92-93 1981/82 59-73

\art M2 B. Malgrange: Introduction aux travaux de J. \'Ecalle! Ens.
Math.! 31 1985 261-282

\art S A.A. Shcherbakov: Topological classification of germs of conformal
mappings with identity linear part! Moscow Univ. Math. Bull.! 37 1982
60-65

\art V S.M. Voronin: Analytic classification of germs of conformal maps
$(\C,0)\to (\C,0)$ with identity linear part! Func. Anal. Appl.! 15 1981
1-17

\bye